\documentclass[11pt,a4paper]{amsart}

\usepackage[utf8]{inputenc}
\usepackage{amsmath}
\usepackage{amsfonts}
\usepackage{amssymb}
\usepackage{amsthm}
\usepackage{setspace}
\usepackage{hyperref}
\usepackage{cleveref}
\usepackage{graphicx}
\usepackage{stmaryrd}
\usepackage[shortlabels]{enumitem}
\usepackage{mathtools}
\usepackage{subcaption}
\usepackage[T1]{fontenc}
\usepackage{xcolor}
\usepackage{float}
\usepackage{tikz-cd}
\usepackage{tikz}
\usepackage[noadjust]{cite}
\usepackage{cite}
\usepackage{ragged2e}
\usepackage[margin=1.2in]{geometry}
\graphicspath{ {./images/} }
\newtheorem{theorem}{Theorem}[section]
\newtheorem{proposition}[theorem]{Proposition}

\theoremstyle{definition}

\newtheorem{remark}[theorem]{Remark}

\newcommand{\ZZ}{\mathbb{Z}}

\newcommand{\bb}[2][]{\ensuremath{\mathbb{#2}^{#1}}}

\newcommand{\GL}{\mathrm{GL}}

\subjclass{20F55, 46L10, 20F65, 46L05}
\keywords{Coxeter groups, von Neumann algebras, strong solidity}

\title{Strong solidity classification of Coxeter groups}

\author[M.~Blufstein]{Mart\'in Blufstein}
\address[M.~Blufstein]{Departamento de Matemática-IMAS (CONICET),
FCEyN,
Universidad de Buenos Aires,
Buenos Aires,
Argentina
}
    
\email{mblufstein@dm.uba.ar}

\author[K.~Goldman]{Katherine Goldman
}
    
    \address[K.~Goldman]{
    Department of Mathematics and Statistics,
    McGill University,
    Burnside Hall,
    805 Sherbrooke Street West,
    Montreal, QC,
    H3A 0B9, Canada}
    
    \email{kat.goldman@mcgill.ca}
    
\author[K.~Oyakawa]{Koichi Oyakawa}

    \address[K.~Oyakawa]{
    Department of Mathematics and Statistics,
    McGill University,
    Burnside Hall,
    805 Sherbrooke Street West,
    Montreal, QC,
    H3A 0B9, Canada}
    
    \email{koichi.oyakawa@mail.mcgill.ca}

\begin{document}

\begin{abstract}
We prove the dichotomy that every Coxeter group either has a strongly solid group von Neumann algebra or contains the product of an infinite cyclic group and a free group of rank 2. This generalizes the same dichotomy for right-angled Coxeter groups by Borst-Caspers. However, our proof is conceptually different, which leads to a significantly streamlined argument. We also provide additional equivalent geometric and group-theoretic characterizations of strong solidity for Coxeter groups that allow us to completely classify those with a strongly solid group von Neumann algebra. In particular, we characterize strong solidity purely in terms of the defining Coxeter-Dynkin diagram. Finally, we obtain the same dichotomy for virtually cocompact special groups.
\end{abstract}

\maketitle 

\section{Introduction}

Strong solidity is a structural property of von Neumann algebras introduced by Ozawa and Popa in \cite{OP10}, where they show that the von Neumann algebras of free groups have this property.
A von Neumann algebra $M$ is called \emph{strongly solid} if for every diffuse amenable von Neumann subalgebra $A \subset M$, the normalizer $\mathcal{N}_M(A)''$, which is generated by the set $\mathcal{N}_M(A)$ of unitaries $u$ of $M$ satisfying $uAu^* = A$, is amenable.  It is known as the strongest indecomposability property that encompasses primeness, solidity, and absence of Cartan subalgebras, and has been extensively studied (see \cite{Hou10,HD11,Tho11,CS13,CSU13,PV14,Iso15,BHV18,Cas21,Cas22,BC24,BCW24}).

Arguably the simplest obstruction for a group to have a strongly solid group von Neumann algebra is to contain $\ZZ \times F_2$, the product of an infinite cyclic group and a free group of rank 2, as a subgroup.
In \cite{BC24}, Borst-Caspers proved a dichotomy that this is the only obstruction in the class of right-angled Coxeter groups: namely, every right-angled Coxeter group either has a strongly solid group von Neumann algebra or contains $\ZZ \times F_2$.
They proved it using various von Neumann algebraic techniques, including Popa's intertwining-by-bimodules theory \cite{popa2006a,popa2006b} and rigidity of amalgamated free products by Vaes \cite{Vae14} and Ioana \cite{Ioa15}. Their method was later applied to obtain more rigidity results in the more general settings of graph products in \cite{BCC25}.

Right-angled Coxeter groups form a (particularly nice) subclass of the general Coxeter groups. 
Coxeter groups arise naturally in geometry as groups generated by reflections with respect to some (not necessarily positive definite) bilinear form on $\bb[n]{R}$. 
They are ubiquitous in mathematics, and play a fundamental role in, for example, the classification of regular polytopes, the classification of semi-simple Lie algebras, the definition of Kazhdan–Lusztig polynomials, the study of braid and Artin groups, representation theory, the geometry of Riemannian symmetric spaces, the topology of aspherical manifolds, and more.
They are also prevalent in other areas of science such as physics, biology, chemistry and computer science.

In this paper, we generalize and strengthen the Borst-Caspers' dichotomy to the class of all Coxeter groups from a new perspective.
Instead of von Neumann algebraic techniques, we use geometric and analytic properties of groups, which results in a very straightforward proof.
The following is our main theorem.

\begin{theorem}\label{thm:main}
    Let $(W,S)$ be a Coxeter system with $S$ finite. The following are equivalent.
    \begin{enumerate}
        \item The group von Neumann algebra of $W$ is strongly solid.
        \item $W$ is biexact.
        \item $W$ does not contain $\ZZ \times F_2$ as a subgroup.
        \item $S$ does not contain two subsets $J_1,J_2$ such that $[J_1,J_2] = 1$, $|\langle J_1 \rangle| = \infty$, and $\langle J_2 \rangle$ is non-elementary word hyperbolic (where $[J_1,J_2]=\{\,sts^{-1}t^{-1} : s \in J_1,\,t \in J_2\,\}$).
        \item $W$ is hyperbolic relative to a finite (possibly trivial) collection of virtually abelian subgroups.
        \item $W$ is hyperbolic relative to a finite (possibly trivial) collection of amenable subgroups.
    \end{enumerate}
\end{theorem}
Recall that a group is ``virtually $\mathcal P$'' for some property $\mathcal P$ if it has a finite index subgroup which is $\mathcal P$.

\begin{remark}
    Coxeter-Dynkin diagrams of finite (\cite[Table 6.1]{Dav25}) and non-elementary hyperbolic (\cite[Theorem 17.1]{moussong1988}) Coxeter groups are classified.
    Hence, by equivalence (4) it is easy to check whether a Coxeter group satisfies any of the equivalences in \Cref{thm:main} by looking at its diagram (see \Cref{sec:prelim} for definitions).
\end{remark}

Biexactness is a functional-analytic generalization of amenability. So the above theorem can be seen as a generalization of the classification of amenable Coxeter groups (see \Cref{prop:amen Coxeter}). In particular, the equivalence of (2) and (3) generalizes the fact that a Coxeter group is amenable if and only if it does not contain $F_2$ as a subgroup.
The equivalence of (3) and (5), while not entirely a new result \cite{Cap09}, can also be seen as generalizing the classification of word hyperbolic Coxeter groups; namely, a Coxeter group is word hyperbolic if and only if it contains no $\bb{Z} \times \bb{Z}$ subgroup \cite[Cor.~12.6.3]{Dav25}.

Our approach relies on geometric properties that are not unique to Coxeter groups.
Consequently, we can extend the classification to the class of virtually cocompact special groups, introduced by Haglund and Wise in \cite{HaglundWise2008}.
A group is called \emph{special} if it is the fundamental group of a special non-positively curved cube complex.
This class plays a central role in geometric group theory, particularly in the study of 3-manifolds and subgroup separability.
Examples of virtually cocompact special groups include word hyperbolic Coxeter groups and right-angled Coxeter groups \cite{NR03}, right-angled Artin groups \cite[Ex.3.3]{HaglundWise2008}, and graph products with finite or cyclic vertex groups \cite[Thm.~3]{Kim12}.
In particular, the following theorem also recovers Borst and Caspers' dichotomy.

\begin{theorem}\label{thm:special}
    Let $G$ be a virtually cocompact special group. The following are equivalent.
    \begin{enumerate}
        \item The group von Neumann algebra of $G$ is strongly solid.
        \item $G$ is biexact.
        \item $G$ does not contain $\ZZ \times F_2$ as a subgroup.
        \item $W$ is hyperbolic relative to a finite (possibly trivial) collection of virtually abelian subgroups.
        \item $W$ is hyperbolic relative to a finite (possibly trivial) collection of amenable subgroups.
    \end{enumerate}
\end{theorem}

It is an open problem whether there exists a non-biexact group whose group von Neumann algebra is solid \cite{DP23}.
\Cref{thm:main} and \ref{thm:special} also show there is no such group in the classes of Coxeter groups or virtually cocompact special groups (see \cite{Oza04}).

\subsection*{Acknowledgements}
We thank Anthony Genevois for bringing \Cref{prop:genevois} to our attention. The first author is a researcher of CONICET. The second author is partially supported by NSF grant DMS-2402105.

\section{Preliminaries} \label{sec:prelim}

In this section, we will primarily introduce basic notation and properties about Coxeter groups. For details on biexact groups, see \cite[Section 15.2]{BO08}, for details on relatively hyperbolic groups, see \cite{Osi06}, and for details on special groups see \cite{HaglundWise2008}.

A Coxeter group can be defined in many equivalent ways, but the way that shall be most convenient for us is via \emph{Coxeter-Dynkin diagrams}.
These are a generalization of Dynkin diagrams used, e.g., to define Weyl groups.
Let $\Gamma$ be a simplicial graph on a finite vertex set $S$; that is, $\Gamma$ is an unoriented graph with no double edges and no loops (edges which start and end at the same vertex).
An edge of $\Gamma$ can then be uniquely represented by an unordered pair $\{s,t\}$ of vertices $s,t \in S$. 
To each edge $\{s,t\}$ we assign a label $m_{st} \in \bb{Z}_{\geq 3} \cup \{ \infty \}$. 
Sometimes if $m_{st} = 3$, we omit the label on the edge when drawing the diagram.
If $s$ and $t$ do not span an edge, then we set $m_{st} = 2$.
For all $s \in S$, we let $m_{ss} = 1$. 
This data gives a Coxeter group $W$ defined by
\[
W = \langle\, S \mid (st)^{m_{st}} = 1 \text{ for all } s, t \in S \text{ with } m_{st} < \infty \,\rangle.
\]
Note that each generator $s \in S$ has order 2, since $s^2 = (ss)^{m_{ss}} = 1$. 
The choice of generators is significant, so we usually notate the choice of generators $S$ along with the group $W$, and call $(W,S)$ a \emph{Coxeter system}.
The Coxeter group $W$ and choice of generators $S$ is actually enough to recover the information encoded in the Coxeter-Dynkin diagram $\Gamma$ (\cite[Ch. IV, \S 1]{bourbaki2002lie}).
We reiterate for emphasis that in our definition, and in the rest of the article, $S$ is always a finite set. 

Suppose $T \subseteq S$. Then we call $W_T \coloneqq \langle T \rangle$ a \emph{(standard) parabolic subgroup} of $W$.
Let $\Gamma(T)$ be the full (or ``induced'') subgraph of $\Gamma$ on vertices $T$, inheriting the labels $m_{st}$ of $\Gamma$.
Let $W'$ be the Coxeter group defined by $\Gamma(T)$.
Then it is a standard fact (see, e.g., \cite[Ch. IV, \S 1, Theorem 2]{bourbaki2002lie}) that $(W', T)$ and $(W_T, T)$ are isomorphic \emph{as Coxeter systems} ($W'$ and $W_T$ are isomorphic via a map which restricts to the identity on $T$).
We say that $(W,S)$ is \emph{irreducible} if there do not exist disjoint $T_1, T_2 \subseteq S$ such that $W = W_{T_1} \times W_{T_2}$.
This is equivalent to requiring that $\Gamma$ be connected.

Coxeter groups have a natural representation as discrete subgroups of $\GL_n(\bb{R})$ generated by (not necessarily orthogonal) reflections, called the \emph{Tits representation} (\cite[Ch. V, \S 4]{bourbaki2002lie}), defined as follows.
For each $s \not= t$ in $S$, define $c_{st} = - \cos(\pi/m_{st})$, and define $c_{ss} = 1$.
Then let $C = (c_{st})_{s,t \in S}$.
This is called the \emph{cosine matrix} of the Coxeter system $(W,S)$. Note that $C$ is symmetric.
Choose a basis $\{e_s\}_{s \in S}$ for \bb[|S|]{R}, and let $B$ be the (symmetric) bilinear form associated to $C$ in this basis, i.e., $B(e_s,e_t) = c_{st}$. 
Then for $s \in S$, we define $\rho_s : (x \mapsto x - 2B(x,e_s)e_s)$.
Then the map $\rho : S \to \GL_n(\bb{R})$ given by $s \mapsto \rho_s$ extends to an injective homomorphism $\rho : W \to \GL_n(\bb{R})$.

The geometry of the bilinear form $B$ has a very close relationship with the  structure of $(W,S)$.
Let $n = |S|$.
The following can be found, e.g., in \cite[\S 6]{Dav25}.
\begin{itemize}
    \item $B$ is positive definite if and only if $|W| < \infty$.
    In this case, $\rho(W)$ is a finite subgroup of $O(n)$, and in particular acts as a (finite) reflection group on the sphere \bb[n-1]{S}.
    For this reason, if $|W| < \infty$, we call $W$ \emph{spherical}.
    \item $B$ is corank-1 if and only if $W$ (setwise) preserves a codimension-1 affine subspace of \bb[n]{R}.
    In this case, it acts as a (cocompact) Euclidean reflection group on this affine hyperplane.
    For this reason, we call such $W$ \emph{affine}.
    \item $B$ is type $(n-1,1)$ if and only if it preserves an (isometrically embedded) copy of \bb[n-1]{H} in \bb[n-1,1]{R}.
    In this case, it acts as a (cocompact) hyperbolic reflection group on this copy of \bb[n-1]{H}.
    For this reason, we call such $W$ \emph{hyperbolic}\footnote{We note that there are many inequivalent notions of ``hyperbolic'' Coxeter group. In particular, there are Coxeter groups which are word hyperbolic (in the sense of Gromov) but are not necessarily hyperbolic in this sense. Although, our notion of hyperbolic Coxeter group does imply word hyperbolicity (and in particular, non-elementary word hyperbolicity). We will always make the distinction clear when it is relevant.}.
\end{itemize}
In any of these cases, we may refer to $S$ or $\Gamma$ themselves as \emph{spherical}, \emph{affine}, or \emph{hyperbolic}, respectively.

While the spherical Coxeter groups have an obvious algebraic characterization (that is, being finite), the affine and hyperbolic Coxeter groups have less obvious consequences.
For the affine Coxeter groups, this is clarified by the following proposition.
First, one more definition: if there exist pairwise disjoint $T_1,\dots,T_n \subseteq S$ such that $W = W_{T_1} \times \dots \times W_{T_n}$ and each $W_{T_i}$ is either spherical or affine, we call $W$ \emph{Euclidean}.
As before, we may sometimes call $S$ or $\Gamma$ themselves \emph{Euclidean}.
\begin{proposition}\label{prop:amen Coxeter}
\textup{\cite[Prop.~17.2.1]{Dav25}}
    Let $(W,S)$ be a Coxeter system. The following are equivalent.
    \begin{enumerate}
        \item $W$ is Euclidean.
        \item $W$ is virtually abelian.
        \item $W$ does not contain a free group of rank $2$.
        \item $W$ is amenable.
    \end{enumerate}
\end{proposition}

\begin{remark}\label{rem:infinite_order_element}
Any virtually solvable subgroup of $W$ is finitely generated (see \cite[Theorem 12.3.4]{Dav25}).
Thus, the equivalence of (2) and (3) implies that every infinite Coxeter group contains an element of infinite order.    
\end{remark}

Suppose $(W,S)$ is not spherical and not affine, but for every maximal proper $T \subset S$ (i.e. $|T| = |S|-1$), $W_T$ is either spherical or irreducible affine.
Then we call $W$ \emph{minimal hyperbolic}. Like with spherical and affine, sometimes we call $T$ or $\Gamma$ themselves \emph{minimal hyperbolic}.
It is worth to note that such $W$ are hyperbolic (and hence word hyperbolic), but not all hyperbolic Coxeter groups are minimal hyperbolic.
Although, every infinite non-affine Coxeter group $(W,S)$ does contain a minimal hyperbolic $W_T$ (by taking $T$ to be a minimal non-affine subset of $S$). 

There is a complete classification of the irreducible spherical, affine, and minimal hyperbolic Coxeter systems in terms of their Coxeter-Dynkin diagrams (\cite[Tables 6.1 \& 6.2]{Dav25}).

For some of the equivalences, it becomes useful to introduce the following notation. If $T \subseteq S$, then define
\[
    T^\perp = \{\, s \in S \setminus T : sts^{-1}t^{-1} = 1 \text{ for all } t \in T \,\}.
\]
The following result of Caprace characterizes when Coxeter groups are relatively hyperbolic with respect to virtually abelian subgroups which are not virtually cyclic.

\begin{proposition} \textup{\cite[Cor.~D]{Cap09}} \label{prop:caprace}
    For a Coxeter system $(W,S)$, the following assertions are equivalent:
  \begin{itemize}
    \item[\textup{(ii)}] For each minimal hyperbolic $T \subset S$, the set $T^\perp$ is spherical.
    \item[\textup{(iv)}] The group $W$ is relatively hyperbolic with respect to a collection of virtually abelian
    subgroups of rank at least~$2$.
  \end{itemize}
\end{proposition}
We note the errata \cite{Cap15} to the article \cite{Cap09}, but this result remains unchanged.

A similar characterization of relative hyperbolicity for virtually cocompact special groups was obtained by Genevois .

\begin{proposition} \textup{\cite[Thm.~1.4]{Gen21}} \label{prop:genevois}
    A virtually cocompact special group is hyperbolic relative to a finite collection of virtually abelian subgroups if and only if it does not contain $\mathbb{Z} \times F_2$ as a subgroup.
\end{proposition}

\section{Proofs of the Theorems}

Using \Cref{rem:infinite_order_element} and \Cref{prop:caprace}, the Theorem follows quickly.

\begin{proof}[Proof of \Cref{thm:main}]
\
\linebreak
\vspace{-15pt}
    \begin{enumerate}[leftmargin=4.5\parindent]
        \item[(3) $\Rightarrow$ (4):] Suppose to the contrapositive that there exist $J_1, J_2 \subseteq S$ such that $[J_1, J_2] = 1$, with $W_{J_1}$ infinite and $W_{J_2}$ non-elementary word hyperbolic.
        Then in particular $[W_{J_1}, W_{J_2}] = 1$.
        Since $|W_{J_1}| = \infty$, by \Cref{rem:infinite_order_element} $W_{J_1}$ contains an element of infinite order, i.e., contains a copy of $\bb{Z}$. Since $W_{J_2}$ is non-elementary word hyperbolic, $W_{J_2}$ contains a copy of the free group $F_2$ of rank 2 (by a ``ping-pong'' argument).
        Then the copy of \bb{Z} in $W_{J_1}$ commutes with the copy of $F_2$ in $W_{J_2}$, so their intersection is trivial as it is central in $F_2$.
        Therefore, together they form a subgroup isomorphic to $\bb{Z} \times F_2$ sitting inside of $W_{J_1 \cup J_2}$.
        
        \item[(4) $\Rightarrow$ (5):] We also show this by contrapositive. Suppose (5) does not hold. Then in particular, condition (iv) of \Cref{prop:caprace} does not hold. By the equivalence of (ii) and (iv) of \Cref{prop:caprace}, there exists a minimal hyperbolic $J \subset S$ such that $J^\perp$ is non-spherical. 
        Then $W_J$ is non-elementary word hyperbolic, $|W_{J^\perp}| = \infty$, and $[J, J^\perp] = 1$ (by the definition of $J^\perp$).
        
        \item[(5) $\Rightarrow$ (6):] This follows since every virtually abelian group is amenable.
        
        \item[(6) $\Rightarrow$ (2):] This follows by \cite[Theorem 1.1]{Oyakawabiexact} (see also \cite{Oza06}), since every amenable group is biexact. 
        
        \item[(2) $\Rightarrow$ (1):] Note that every Coxeter group is weakly amenable.
        This is because every Coxeter group admits a proper action on a locally finite, finite dimensional CAT(0) cube complex \cite[Thm.~1]{NR03}, and every such group is weakly amenable by \cite{GN10}.
        By \cite[Cor.~0.2]{CSU13}, the group von Neumann algebra of a biexact weakly amenable group is strongly solid.
        
        \item[(1) $\Rightarrow$ (3):] Finally, this follows since strong solidity passes down to von Neumann subalgebras, and the group von Neumann algebra of $\ZZ \times F_2$ is not strongly solid.
        For more detail, one can argue as follows. For any group $G$, let $\mathcal{L}G$ denote the group von Neumann algebra of $G$.
        If $W$ contains $\ZZ\times F_2$, then the normalizer $\mathcal{N}_{\mathcal{L}W}(\mathcal L \ZZ)''$ contains $\mathcal{L}F_2$, because $\mathcal{L}\ZZ$ and $\mathcal{L}F_2$ commute.
        Hence, $\mathcal{L}W$ is not strongly solid since $\mathcal{L}\ZZ$ is diffuse amenable and $\mathcal{L}F_2$ is non-amenable, so $\mathcal{N}_{\mathcal{L}W}(\mathcal L \ZZ)''$ is non-amenable. \qedhere
    \end{enumerate}
\end{proof}

The proof of \Cref{thm:special} is exactly the same, the only difference being that \Cref{prop:genevois} should be used in place of \Cref{prop:caprace}.

\bibliographystyle{amsalpha}
\bibliography{mybib}

\end{document}